\documentclass{article}
\usepackage{graphicx}
\usepackage{amsmath}
\usepackage{amsfonts}
\usepackage{amssymb}
\newtheorem{theorem}{Theorem}

\newtheorem{corollary}[theorem]{Corollary}

\newtheorem{definition}[theorem]{Definition}

\newtheorem{lemma}[theorem]{Lemma}
\newtheorem{notation}[theorem]{Notation}

\newtheorem{proposition}[theorem]{Proposition}
\newtheorem{remark}[theorem]{Remark}

\newenvironment{proof}[1][Proof]{\textbf{#1.} }{\ \rule{0.5em}{0.5em}}

\begin{document}

\title{On $k$-abelian $p$-filiform Lie algebras}
\author{Otto Rutwig Campoamor Stursberg\\Departamento de Geometr\'{\i}a y Topolog\'{\i}a\\Fac. CC. Matem\'{a}ticas\\Universidad Complutense\\28040 Madrid ( Spain )\\rutwig@sunal1.mat.ucm.es}
\date{}
\maketitle
\begin{abstract}
We classify the $\left(  n-5\right)  $-filiform Lie algebras which have the
additional property of a non-abelian derived subalgebra. We show that this
property is strongly related with the structure of the Lie algebra of
derivations; explicitely we show that if a $\left(  n-5\right)  $-filiform Lie
algebra is characteristically nilpotent, then it must be $2$-abelian. We also
give applications of $k$-abelian Lie algebras to the construction of solvable
rigid laws whose nilradical is $k$-abelian with mixed characteristic sequence,
as well as applications to the theory of nilalgebras of parabolic subalgebras
of the exceptional simple Lie algebra $E_{6}$.
\end{abstract}

\section{Generalities}

\begin{definition}
Let $\frak{g}$ be a finite dimensional vectorial space over $\mathbb{C}$. A
Lie algebra law over $\mathbb{C}^{n}$ is a bilinear alternated mapping $\mu\in
Hom\left(  \mathbb{C}^{n}\times\mathbb{C}^{n},\mathbb{C}^{n}\right)  $ which
satisfies the conditions

\begin{enumerate}
\item $\mu\left(  X,X\right)  =0,\;\forall\;X\in\mathbb{C}^{n}$

\item $\mu\left(  X,\mu\left(  Y,Z\right)  \right)  +\mu\left(  Z,\mu\left(
X,Y\right)  \right)  +\mu\left(  Y,\mu\left(  Z,X\right)  \right)
=0,\;\forall\;X,Y,Z\in\mathbb{C}^{n}$, \ \newline ( Jacobi identity )
\end{enumerate}

If $\mu$ is a Lie algebra law, the pair $\frak{g}=\left(  \mathbb{C}^{n}%
,\mu\right)  $ is called Lie algebra. From now on we identify the Lie algebra
with its law $\mu$.
\end{definition}

\begin{remark}
We say that $\mu$ is the law of $\frak{g}$, and where necessary we use the
bracket notation to describe the law :
\[
\left[  X,Y\right]  =\mu\left(  X,Y\right)  ,\;\forall\;X,Y\in\frak{g}%
\]
The nondefined brackets are zero or obtained by antisymmetry.
\end{remark}

Let $\frak{g}_{n}=\left(  \mathbb{C}^{n},\mu\right)  $ be a nilpotent Lie
algebra. For any nonzero vector $X\in\frak{g}_{n}-C^{1}\frak{g}_{n}$ let
$c\left(  X\right)  $ be the ordered sequence of a similitude invariants for
the nilpotent operator $ad_{\mu}\left(  X\right)  \dot{,}$ i.e., the ordered
sequence of dimensions of Jordan blocks of this operator. The set of these
sequences is ordered lexicographically.

\begin{definition}
The characteristic sequence of $\frak{g}_{n}$ is an isomorphism invariant
$c\left(  \frak{g}_{n}\right)  $ defined by
\[
c\left(  \frak{g}_{n}\right)  =\max_{X\in\frak{g}_{n}-C^{1}\frak{g}_{n}%
}\left\{  c\left(  X\right)  \right\}
\]
A nonzero vector $X\in\frak{g}_{n}-C^{1}\frak{g}_{n}$ for which $c\left(
X\right)  =c\left(  \frak{g}_{n}\right)  $ is called characteristic vector.
\end{definition}

\begin{remark}
In particular, the algebras with maximal characteristic sequence $\left(
n-1,1\right)  $ correspond to the filiform algebras introduced by Vergne
[Ve1]. Thus it is natural to generalize this concept to lower sequences.
\end{remark}

\begin{definition}
A nilpotent Lie algebra $\frak{g}_{n}$ is called $p$-filiform if its
characteristic sequence is $\left(  n-p,1,..^{p}..,1\right)  .$
\end{definition}

\begin{remark}
This definition was first given in [CGoJ]. The $\left(  n-1\right)  $-filiform
Lie algebras are abelian, while the $\left(  n-2\right)  $-filiform are the
direct sum of an Heisenberg algebra $\frak{h}_{2p+1}$ and abelian algebras. A
classification of the $\left(  n-3\right)  $-filiform can also be found in
[CGoJ]. These and any $\left(  n-4\right)  $-filiform Lie algebra have non
trivial diagonalizable derivations [AC1]. This fact is important, for it is
telling us that their structure is relatively simple. To search for nilpotent
algebras with rank zero ( i.e, with no nonzero diagonalizable derivations ) we
must start with the $\left(  n-5\right)  $-filiform Lie algebras. As the
difficulty of distinguishing isomorphism classes increases considerably for
bigger indexes, it seems reasonable to consider additional assumptions made on
the algebras to be classified. For example, the filtration given by the
central descending sequence can be used to impose additional conditions on the
$p$-filiformness.
\end{remark}

\begin{remark}
For indexes $p\geq\left(  n-4\right)  $ the number of isomorphism classes is
finite. The index $p=\left(  n-5\right)  $ is the first for which an infinity
of isomorphism classes exists.
\end{remark}

\begin{definition}
Let $\frak{g}_{n}$ be a nilpotent Lie algebra. The smallest integer $k$ such
that the ideal $C^{k}\frak{g}_{n}$ is abelian is called commutativity index of
$\frak{g}_{n}.$
\end{definition}

\begin{definition}
A Lie nilpotent algebra $\frak{g}_{n}$ is called $k$-abelian if $k$ is the
smallest positive integer such that
\[
C_{\frak{g}}\left(  C^{k}\frak{g}\right)  \supset C^{k}\frak{g}\text{\ \ \ and
\ \ }C^{k-1}\frak{g}\supsetneq C^{k-1}\frak{g}%
\]
\end{definition}

\begin{remark}
The preceding definition is equivalent to impose that the commutativity index
of $\frak{g}_{n}$ is exactly $k$. In [GGoK] a less restrictive definition of
$k$-abelianity is considered. The purpose there is to study certain
topological properties of the variety of filiform laws $\frak{F}^{m}.$ Our
definition is more restrictive: the $k$-abelian Lie algebras do \textbf{not}
contain the $\left(  k-1\right)  $-abelian algebras; the reason is justified
by the important structural difference between algebras having its ideal
$C^{k}\frak{g}_{n}$ abelian and those having it not. On the other side we
avoid unnecessary repetitions.
\end{remark}

As we are considering here the $\left(  n-5\right)  $-filiform Lie algebras,
we have to determine which abelianity indexes are admissible. Only the
nonsplit algebras are of interest for us, thus from now on we will understand
nonsplit Lie algebra \ when we say Lie algebra.

\begin{lemma}
Let $\frak{g}_{n}$ be an $\left(  n-5\right)  $-filiform Lie algebra. Then
there exists a basis $\left\{  \omega_{1},..,\omega_{6},\theta_{1}%
,..,\theta_{n-6}\right\}  $ of $\left(  \mathbb{C}^{n}\right)  ^{\ast}$ \ such
that the law is expressible as
\begin{align*}
d\omega_{1} &  =d\omega_{2}=0\\
d\omega_{3} &  =\omega_{1}\wedge\omega_{2}\\
d\omega_{4} &  =\omega_{1}\wedge\omega_{3}+\sum_{i=1}^{n-6}\alpha_{i}%
^{2}\theta_{i}\wedge\omega_{2}\\
d\omega_{5} &  =\omega_{1}\wedge\omega_{4}+\sum_{i=1}^{n-6}\left(  \alpha
_{i}^{2}\theta_{i}\wedge\omega_{3}+\alpha_{i}^{3}\theta_{i}\wedge\omega
_{2}\right)  +\beta_{2}\omega_{3}\wedge\omega_{2}\\
d\omega_{6} &  =\omega_{1}\wedge\omega_{5}+\sum_{i=1}^{n-6}\left(  \alpha
_{i}^{2}\theta_{i}\wedge\omega_{4}+\alpha_{i}^{3}\theta_{i}\wedge\omega
_{3}+\alpha_{i}^{4}\theta_{i}\wedge\omega_{2}\right)  \\
&  +\sum_{1\leq i,j\leq n-6}a_{ij}^{1}\theta_{i}\wedge\theta_{j}+\beta
_{1}\left(  \omega_{5}\wedge\omega_{2}-\omega_{3}\wedge\omega_{4}\right)
+\beta_{2}\omega_{4}\wedge\omega_{2}+\beta_{3}\omega_{3}\wedge\omega_{2}\\
d\theta_{j} &  =\varepsilon_{i}^{j}\theta_{i}\wedge\omega_{2}-\beta
_{1,j}\left(  \omega_{5}\wedge\omega_{2}-\omega_{3}\wedge\omega_{4}\right)
+\beta_{3,j}\omega_{3}\wedge\omega_{2},\;1\leq j\leq n-6
\end{align*}
\end{lemma}

The proof is trivial.

\begin{lemma}
Any $\left(  n-5\right)  $-filiform Lie algebra $\frak{g}_{n}$ is either $1$
or $2$-abelian.
\end{lemma}

\begin{proof}
If the algebra is $1$-abelian, it is simply an algebra whose derived algebra
is abelian [Bra]. If it is $2$-abelian, then there exist $X,Y\in C^{1}%
\frak{g}_{n}$ such that $0\neq\lbrack X,Y]$. From the above equations it is
immediate to derive the possibilities:

\begin{enumerate}
\item $\dim\;C^{1}\frak{g}_{n}=6$ \ 

\item $\dim\;C^{1}\frak{g}_{n}=5$ and $\exists X,Y\in C^{1}\frak{g}_{n}$ such
that $[X,Y]\neq0$

\item $\dim\;C^{1}\frak{g}_{n}=4$ and $\exists X,Y\in C^{1}\frak{g}_{n}$ such
that $[X,Y]\neq0.$\newline 
\end{enumerate}
\end{proof}

\begin{remark}
From this lemma we see how important is to consider our stronger version of
the $k$-abelianity. In particular we will see its connection with the
characteristic nilpotence.
\end{remark}

\begin{notation}
For $n\geq7$ let $\frak{g}_{0}^{n}$ be the Lie algebra whose Cartan-Maurer
equations are
\begin{align*}
d\omega_{1}  &  =d\omega_{2}=0\\
d\omega_{j}  &  =\omega_{1}\wedge\omega_{j-1},\;3\leq j\leq6\\
d\theta_{j}  &  =0,\;1\leq j\leq n-6
\end{align*}
\end{notation}

Let $\{X_{1},..,X_{6},Y_{1},..,Y_{n-6}\}$ be a dual basis of $\{\omega
_{1},..,\omega_{6},\theta_{1},..,\theta_{6}\}$. Let $V_{1}=\left\langle
X_{1},..,X_{6}\right\rangle _{\mathbb{C}}$ and $V_{2}=\left\langle
Y_{1},..,Y_{n-6}\right\rangle _{\mathbb{C}}$. We write $B\left(  V_{i}%
,V_{j}\right)  $ to denote the space of bilinear alternated mappings from
$V_{i}$ to $V_{j}$.\newline Let us consider the following applications :

\begin{enumerate}
\item $\psi_{i,j}^{1}\in B\left(  V_{2},V_{1}\right)  ,\;1\leq i,j\leq n-6$ :
\[
\psi_{i,j}^{1}\left(  Y_{i},Y_{l}\right)  =\left\{
\begin{array}
[c]{cc}%
X_{6} & \text{if \ }i=k\;,j=l\\
0 & \text{otherwise}%
\end{array}
\right.
\]

\item $\psi_{i}^{j}\in Hom\left(  V_{2}\times V_{1},V_{1}\right)  ,\;j=2,3,4$%
\[
\psi_{i}^{j}\left(  Y_{k},X_{l}\right)  =\left\{
\begin{array}
[c]{cc}%
X_{l+j} & \text{if \ }i=k,\;2\leq l\leq6-j\\
0 & \text{otherwise}%
\end{array}
\right.
\]

\item $\varphi_{1,k}\in B\left(  V_{1},V_{2}\right)  ,\;1\leq k\leq n-6$ :%
\[
\varphi_{1,k}\left(  X_{5},X_{2}\right)  =\varphi_{1,k}\left(  X_{3}%
,X_{4}\right)  =Y_{k}%
\]

\item $\varphi_{3,k}\in B\left(  V_{1},V_{2}\right)  ,\;1\leq k\leq n-6$ :%
\[
\varphi_{3,k}\left(  X_{3},X_{2}\right)  =Y_{k}%
\]

\item $\varphi_{1}\in B\left(  V_{1},V_{1}\right)  $ :%
\[
\varphi_{1}\left(  X_{5},X_{2}\right)  =\varphi_{1}\left(  X_{3},X_{4}\right)
=X_{6}%
\]

\item $\varphi_{2}\in B\left(  V_{1},V_{1}\right)  $ :%
\[
\varphi_{\uparrow}\left(  X_{k},X_{2}\right)  =X_{k+2},\;k=3,4
\]

\item $\varphi_{3}\in B\left(  V_{1},V_{1}\right)  $ :%
\[
\varphi_{3}\left(  X_{3},X_{2}\right)  =X_{6}%
\]
\end{enumerate}

where the undefined brackets are zero or obtained by antisymmetry.

\begin{lemma}
For $n\geq7$,\ $1\leq k\leq n-6$ and $l=2,3,4$ the mappings \ $\psi_{1,k}%
,\psi_{3,k},\psi_{i}^{l},$\newline $\varphi_{1,k},\varphi_{3,k},\varphi
_{1},\varphi_{2},\varphi_{3}$ are $2$-cocycles of the subspace $Z^{2}\left(
\frak{g}_{0}^{n},\frak{g}_{0}^{n}\right)  $.
\end{lemma}

\begin{notation}
For convenience in the exposition, we introduce the following notation
\[
\sum_{\substack{t=i\\m>k}}^{j}\psi_{f\left(  t\right)  ,f\left(  t\right)  +1}%
\]
where the sum is only defined whenever $m\geq k+1$.
\end{notation}

\begin{proposition}
Any nonsplit $\left(  n-5\right)  $-filiform Lie algebra with $\dim
\,C^{1}\frak{g}_{n}=6$ is isomorphic to one of the following laws :\newline 

\begin{enumerate}
\item $\frak{g}_{2m}^{1}\left(  m\geq4\right)  :$%
\[
\frak{g}_{0}^{2m}+\varphi_{1,1}+\varphi_{3,2}+\psi_{2}^{3}+\sum
_{\substack{t=2\\m>4}}^{m-3}\psi_{2t-1,2t}^{1}%
\]
\newline 

\item $\frak{g}_{2m}^{2}\left(  m\geq5\right)  :$%
\[
\frak{g}_{0}^{2m}+\varphi_{1,1}+\varphi_{3,2}+\psi_{3}^{3}+\sum_{t=2}%
^{m-3}\psi_{2t-1,2t}^{1}%
\]

\item $\frak{g}_{2m+1}^{3}\left(  m\geq4\right)  :$%
\[
\frak{g}_{0}^{2m+1}+\varphi_{1,1}+\varphi_{3,2}+\psi_{3}^{3}+\sum
_{\substack{t=2\\m>4}}^{m-3}\psi_{2t,2t+1}^{1}%
\]

\item $\frak{g}_{2m}^{4}\left(  m\geq4\right)  :$%
\[
\frak{g}_{0}^{2m}+\varphi_{1,1}+\varphi_{3,2}+\sum_{\substack{t=2\\m>4}%
}^{m-3}\psi_{2t-1,2t}^{1}%
\]

\item $\frak{g}_{2m}^{5}\left(  m\geq4\right)  :$%
\[
\frak{g}_{0}^{2m}+\varphi_{1,1}+\varphi_{2}+\varphi_{3,2}+\psi_{2}^{3}%
+\sum_{\substack{t=2\\m>4}}^{m-3}\psi_{2t-1,2t}^{1}%
\]
\end{enumerate}
\end{proposition}

\begin{proof}
Suppose $\beta_{1,k}\neq0$ for $k\neq1$. We can take $\beta_{1,1}=1$ and
$\beta_{1,i}=0$ for $i\geq2$, as well as $\beta_{1}=0$. The Jacobi conditions
imply
\[
\alpha_{1}^{2}=\alpha_{1}^{3}=\alpha_{1}^{1}=0,\;\alpha_{1}^{4}=\sum_{k\geq
2}\beta_{3,k}\alpha_{k,j}^{4},\;\sum\beta_{3,k}\alpha_{i,k}^{1}=0
\]
Let $\beta_{3,k}\neq0$ for $k\neq1$, so that we can choose $\beta
_{3,2}=1,\beta_{3,k}=0$ for $k\neq2$. A change of basis allows to take
$\beta_{3}=0$. From the conditions above we deduce $\alpha_{1}^{4}=\alpha
_{i}^{1}=0$. Consider tha change $\omega_{2}^{\prime}=\alpha\omega_{1}%
+\omega_{2}$ with $\alpha\neq0$. Then we have
\[
\left\{
\begin{array}
[c]{c}%
\;\;\;\alpha_{j,i}^{4}=0,\;j\geq2\\
\alpha_{j}^{3}\beta_{2}=0,\;j\geq2
\end{array}
\right.
\]
There are two cases :

\begin{enumerate}
\item
\begin{enumerate}
\item  If $\alpha_{2}^{3}\neq0$ we suppose $\alpha_{2}^{3}=1$ and $\alpha
_{i}^{3}=0,\;\forall i\neq2$ through a linear change$.$ Reordering the forms
$\theta_{i}$ we can suppose $\alpha_{2t-1,2t}^{1}=1$ for $2\leq t\leq
\frac{n-6}{2};\;\alpha_{i,j}^{1}$ for the remaining. We obtain a unique class
of nonsplit Lie algebras in even dimension and isomorphic to $\frak{g}%
_{2m}^{1}.$

\item  If $\alpha_{2}^{3}=0$

\begin{enumerate}
\item  If \ $\alpha_{i}^{3}\neq0$ with $i\geq3$ we can suppose $\alpha^{3}%
{}_{3}=1$ and $\alpha_{j}^{3}=0$ for $j\neq3.$ Reordering the $\theta_{i}$ we
obtain one algebra in even and one algebra in odd dimension, which are
respectively isomorphic to $\frak{g}_{2m}^{2}$ and $\frak{g}_{2m+1}^{3}.$

\item  If $\alpha_{i}^{3}=0$ for $i\geq3$ we obtain in an analogous way two
even dimensional algebras, respectively isomorphic to $\frak{g}_{2m}^{4}$ and
$\frak{g}_{2m}^{5}.$
\end{enumerate}
\end{enumerate}
\end{enumerate}
\end{proof}

\begin{remark}
It is very easy to see that the obtained algebras are pairwise non isomorphic,
as their infinitesimal deformations are not cohomologous cocycles in the
cohomology space $H^{2}\left(  \frak{g}_{0}^{2m},\frak{g}_{0}^{2m}\right)  $.
This calculations are routine and will be ommited in future.
\end{remark}

\begin{remark}
From the linear system associated ( for the elementary properties of these
systems see [AG1] ) to the algebras above it follows the existence of nonzero
eigenvectors for diagonalizable derivations, so that the rank is at least one.
Then the algebra of derivations has nonzero semi-simple derivations.
\end{remark}

\begin{notation}
We define the set
\[
\frak{h}_{2}=\left\{  \frak{g}\;|\;\frak{g}\text{ \ is nonsplit,
}2\text{-abelian and }\left(  n-5\right)  \text{-filiform }\right\}
\]
\end{notation}

We now express conditions making reference to the reduced system of forms
given before :

We say that $\frak{g}\in\frak{H}_{2}$ satisfies property $\left(  P1\right)  $
if%
\begin{align*}
\dim\;C^{1}\frak{g}  &  =5\\
\beta_{1,1}  &  =1
\end{align*}

\begin{remark}
The general condition would be $\beta_{1,k}\neq0$ for a $k\geq1$ and
$\beta_{3,k}=0$ for any $k$. Now a elementary change of basis allows to reduce
it to the preceding form.
\end{remark}

\begin{proposition}
Let $\frak{g}_{n}$ be an $\left(  n-5\right)  $-filiform Lie algebra
satisfying the property $\left(  P1\right)  .$ The $\frak{g}_{n}$ is
isomorphic to one of the following laws:\newline 

\begin{enumerate}
\item $\frak{g}_{2m+1}^{6}\left(  m\geq4\right)  :$%
\[
\frak{g}_{0}^{2m+1}+\varphi_{1,1}+\psi_{2}^{3}+\sum_{t=2}^{m-3}\psi
_{2t,2t+1}^{1}%
\]

\item $\frak{g}_{2m}^{7}\left(  m\geq4\right)  :$%
\[
\frak{g}_{0}^{2m}+\varphi_{1,1}+\psi_{2}^{3}+\sum_{\substack{t=2\\m>4}%
}^{m-3}\psi_{2t-1,2t}^{1}%
\]

\item $\frak{g}_{2m+1}^{8}\left(  m\geq3\right)  :$%
\[
\frak{g}_{0}^{2m+1}+\varphi_{1,1}+\varphi_{2}+\sum_{\substack{t=2\\m>3}%
}^{m-3}\psi_{2t,2t+1}^{1}%
\]

\item $\frak{g}_{2m+1}^{9}\left(  m\geq3\right)  :$%
\[
\frak{g}_{0}^{2m+1}+\varphi_{1,1}+\varphi_{3}+\sum_{\substack{t=1\\m>3}%
}^{m-3}\psi_{2t,2t+1}%
\]

\item $\frak{g}_{2m+1}^{10}\left(  m\geq3\right)  :$%
\[
\frak{g}_{0}^{2m+1}+\varphi_{1,1}+\sum_{\substack{t=1\\m>3}}^{m-1}%
\psi_{2t,2t+1}^{1}%
\]
\end{enumerate}
\end{proposition}

\begin{proof}
$\beta_{3,k}=0$ for all $k.$ The characteristic sequence implies $\alpha
_{i}^{4}=\alpha_{i}^{3}\beta_{3}=0$ for all $i.$ Moreover, a change of basis
of the type $\omega_{2}^{\prime}=\omega_{4}-\frac{j^{1}}{2}\omega_{2}$ allows
to suppose $\alpha_{1,j}^{4}=0.$

\begin{enumerate}
\item  If $\exists\,\alpha_{i}^{3}\neq0$ we suppose $\alpha_{2}^{3}%
=1,\alpha_{i}^{3}=0,\forall i\neq2.$ A change of basis allows $\beta_{3}=0.$
There are two possibilities: an even dimensional algebra isomorphic to
$\frak{g}_{2m}^{7}$ and an odd dimensional one isomorphic to $\frak{g}%
_{2m+1}^{6}.$

\item $\alpha_{i}^{3}=0,\forall i.$

\begin{enumerate}
\item  If $\beta_{2}\neq0$ we put $\ \beta_{2}=1$ and $\beta_{3}=0$ with a
linear change of basis. We obtain a unique algebra in odd dimension isomorphic
to $\frak{g}_{2m+1}^{8}.$

\item  If $\beta_{2}=0$ there are two possibilities, depending on $\beta_{3}$
: we obtain two odd dimensional algebras which are respectively isomorphic to
$\frak{g}_{2m+1}^{9}$ and $\frak{g}_{2m+1}^{10}.$
\end{enumerate}
\end{enumerate}
\end{proof}

A Lie algebra $\frak{g}\in\frak{h}_{2}$ satisfies property $\left(  P2\right)
$ if
\begin{align*}
\dim C^{1}\frak{g}_{n}  &  =5\\
\beta_{3,t}  &  \neq0\text{ for t}\geq1
\end{align*}

\begin{proposition}
Let $\frak{g}_{n}$ be an algebra with property $\left(  P2\right)  $.
\ Then\ $\frak{g}_{n}$ is isomorphic to one of the following laws\newline 

\begin{enumerate}
\item $\frak{g}_{2m+1}^{11}\left(  m\geq3\right)  :$%
\[
\frak{g}_{0}^{2m+1}+\varphi_{3,1}+\psi_{1}^{3}+\psi_{1}^{4}+\sum
_{\substack{t=1\\m>3}}^{m-3}\psi_{2t,2t+1}^{1}%
\]

\item $\frak{g}_{2m}^{12}\left(  m\geq4\right)  :$%
\[
\frak{g}_{0}^{2m}+\varphi_{3,1}+\psi_{1}^{3}+\psi_{1}^{4}+\psi_{2}^{4}%
+\sum_{\substack{t=2\\m>4}}^{m-3}\psi_{2t-1,2t}^{1}%
\]

\item $\frak{g}_{2m}^{13}\left(  m\geq4\right)  :$%
\[
\frak{g}_{0}^{2m}+\varphi_{3,1}+\psi_{1}^{3}+\psi_{2}^{4}+\sum
_{\substack{t=2\\m>4}}^{m-3}\psi_{2t-1,2t}^{1}%
\]

\item $\frak{g}_{2m+1}^{14}\left(  m\geq4\right)  :$%
\[
\frak{g}_{0}^{2m+1}+\varphi_{3,1}+\psi_{1}^{3}+\sum_{t=1}^{m-3}\psi
_{2t,2t+1}^{1}%
\]

\item $\frak{g}_{2m+1}^{15}\left(  m\geq4\right)  :$%
\[
\frak{g}_{0}^{2m+1}+\varphi_{3,1}+\varphi_{1}+\psi_{1}^{4}+\psi_{2}^{3}%
+\sum_{t=1}^{m-3}\psi_{2t,2t+1}^{1}%
\]

\item $\frak{g}_{2m}^{16}\left(  m\geq4\right)  :$%
\[
\frak{g}_{0}^{2m}+\varphi_{3,1}+\varphi_{1}+\psi_{1}^{4}+\psi_{2}^{3}%
+\sum_{\substack{t=2\\m>4}}^{m-3}\psi_{2t-1,2t}^{1}%
\]

\item $\frak{g}_{2m+1}^{17}\left(  m\geq3\right)  :$%
\[
\frak{g}_{0}^{2m+1}+\varphi_{3,1}+\varphi_{1}+\psi_{1}^{4}+\sum
_{\substack{t=1\\m>3}}^{m-3}\psi_{2t,2t+1}^{1}%
\]

\item $\frak{g}_{2m+1}^{18}\left(  m\geq4\right)  :$%
\[
\frak{g}_{0}^{2m+1}+\varphi_{3,1}+\varphi_{1}+\psi_{3}^{4}+\psi_{2}^{3}%
+\sum_{\substack{t=2\\m>4}}^{m-3}\psi_{2t,2t+1}^{1}%
\]

\item $\frak{g}_{2m}^{19}\left(  m\geq5\right)  :$%
\[
\frak{g}_{0}^{2m}+\varphi_{3,1}+\varphi_{1}+\psi_{24}^{1}+\psi_{2}^{3}%
+\psi_{3}^{4}+\sum_{\substack{t=2\\m>5}}^{m-4}\psi_{2t+1,2t+2}^{1}%
\]

\item $\frak{g}_{2m}^{20}\left(  m\geq4\right)  :$%
\[
\frak{g}_{0}^{2m}+\varphi_{3,1}+\varphi_{1}+\psi_{2}^{3}+\sum
_{\substack{t=2\\m>4}}^{m-3}\psi_{2t-1,2t}^{1}%
\]

\item $\frak{g}_{2m+1}^{21}\left(  m\geq4\right)  :$%
\[
\frak{g}_{0}^{2m+1}+\varphi_{3,1}+\varphi_{1}+\psi_{2}^{3}+\sum_{t=1}%
^{m-3}\psi_{2t,2t+1}^{1}%
\]

\item $\frak{g}_{2m}^{22}\left(  m\geq4\right)  :$%
\[
\frak{g}_{0}^{2m}+\varphi_{3,1}+\varphi_{1}+\psi_{2}^{4}+\sum
_{\substack{t=2\\m>4}}^{m-3}\psi_{2t-1,2t}^{1}%
\]

\item $\frak{g}_{2m+1}^{23}\left(  m\geq3\right)  :$%
\[
\frak{g}_{0}^{2m+1}+\varphi_{3,1}+\varphi_{1}+\sum_{t=1}^{m-3}\psi
_{2t,2t+1}^{1}%
\]

\item $\frak{g}_{2m+1}^{24,\alpha}\left(  m\geq3\right)  :$%
\[
\frak{g}_{0}^{2m+1}+\varphi_{3,1}+\alpha\varphi_{2}+\psi_{1}^{3}+\psi_{1}%
^{4}+\sum_{\substack{t=1\\m>3}}^{m-3}\psi_{2t,2t+1}^{1}%
\]

\item $\frak{g}_{2m}^{25,\alpha}\left(  m\geq4\right)  :$%
\[
\frak{g}_{0}^{2m}+\varphi_{3,1}+\alpha\varphi_{2}+\psi_{1}^{3}+\psi_{1}%
^{4}+\psi_{2}^{4}+\sum_{\substack{t=2\\m>4}}^{m-3}\psi_{2t-1,2t}^{1}%
\]

\item $\frak{g}_{2m}^{26}\left(  m\geq4\right)  :$%
\[
\frak{g}_{0}^{2m}+\varphi_{3,1}+\varphi_{2}+\psi_{1}^{3}+\psi_{2}^{4}%
+\sum_{\substack{t=2\\m>4}}^{m-3}\psi_{2t-1,2t}^{1}%
\]

\item $\frak{g}_{2m+1}^{27}\left(  m\geq3\right)  :$%
\[
\frak{g}_{0}^{2m+1}+\varphi_{3,1}+\varphi_{2}+\psi_{1}^{3}+\sum
_{\substack{t=1\\m>3}}^{m-3}\psi_{2t,2t+1}^{1}%
\]

\item $\frak{g}_{2m+1}^{28}\left(  m\geq4\right)  :$%
\[
\frak{g}_{0}^{2m+1}+\varphi_{3,1}+\varphi_{1}+\varphi_{2}+\psi_{2}^{3}%
+\psi_{3}^{4}+\sum_{\substack{t=2\\m>4}}^{m-3}\psi_{2t,2t+1}^{1}%
\]

\item $\frak{g}_{2m}^{29}\left(  m\geq5\right)  :$%
\[
\frak{g}_{0}^{2m}+\varphi_{3,1}+\varphi_{1}+\varphi_{2}+\psi_{2}^{3}+\psi
_{3}^{4}+\psi_{24}^{1}+\sum_{\substack{t=2\\m>5}}^{m-4}\psi_{2t+1,2t+2}^{1}%
\]

\item $\frak{g}_{2m}^{30}\left(  m\geq4\right)  :$%
\[
\frak{g}_{0}^{2m}+\varphi_{3,1}+\varphi_{1}+\varphi_{2}+\psi_{2}^{3}%
+\sum_{\substack{t=2\\m>4}}^{m-3}\psi_{2t-1,2t}^{1}%
\]

\item $\frak{g}_{2m+1}^{31}\left(  m\geq4\right)  :$%
\[
\frak{g}_{0}^{2m+1}+\varphi_{3,1}+\varphi_{1}+\varphi_{2}+\psi_{2}^{3}%
+\sum_{t=1}^{m-3}\psi_{2t,2t+1}^{1}%
\]

\item $\frak{g}_{2m}^{32}\left(  m\geq4\right)  :$%
\[
\frak{g}_{0}^{2m}+\varphi_{3,1}+\varphi_{1}+\varphi_{2}+\psi_{2}^{4}%
+\sum_{\substack{t=2\\m>4}}^{m-3}\psi_{2t-1,2t}^{1}%
\]
\end{enumerate}
\end{proposition}

\begin{proof}
The starting assumptions are $\beta_{1,k}=0$ for all $k$ and $\beta_{3,k}%
\neq0$ for some $k\geq1.$ We can suppose $\beta_{3,1}=1$, and from the Jacobi
conditions we obtain\
\[
\alpha_{1}^{2}=0,\;\alpha_{i1}^{1}=2\alpha_{i}^{2}\beta_{1},\;i\geq2
\]
From the characteristic sequence we deduce the nullity of $\alpha_{i}^{3}$ for
all $i,$ so that $\alpha_{i1}^{1}=0$ $\forall i$ as well.

\begin{enumerate}
\item $\alpha_{1}^{3}=1$ : a combination of the linear changes of basis allow
to take $\beta_{1}=0.$

\begin{enumerate}
\item $\alpha_{1}^{4}=1$ :.

\begin{enumerate}
\item  If $\alpha_{i}^{4}=0$, $\forall i\geq2$ we reorder the $\left\{
\theta_{2},..,\theta_{n-6}\right\}  $ such that $\alpha_{2t-1,2t}^{1}=1$
$\,$for $1\leq t\leq\frac{n-6}{2}$ and the remaining brackets zero. The
decisive structure constant is $\beta_{2}.$ If it is zero we obtain an odd
dimensional Lie algebra isomorphic to $\frak{g}_{2m+1}^{11}.$ If not,
$\beta_{2}=\alpha$ \ is an essential parameter. So we obtain an infinite
family of odd dimensional Lie algebras isomorphic to the family $\frak{g}%
_{2m+1}^{24,\alpha}.$

\item $\exists\,\alpha_{i}^{4}\neq0$, $i\geq2.$ Without loss of generality we
can choose $\alpha_{2}^{4}\neq0$ and the remaining zero for $i\geq3.$ It is
easy to deduce $\alpha_{2j}^{1}=0,\,\forall j.$ Reordering $\{\theta
_{3},..,\theta_{n-6}\}$ in the previous manner we obtain an even dimensional
Lie algebra and an infinite family of even dimensional algebras, which are
respectively isomorphic to $\frak{g}_{2m}^{12}$ and $\frak{g}_{2m}^{25,\alpha}.$
\end{enumerate}

\item  Now take $\alpha_{1}^{4}=0$

\begin{enumerate}
\item  If there is an index $i\geq2$ such that $\alpha_{i}^{4}\neq0$ we can
suppose $\alpha_{2}^{4}=1$ and the remaining zero. Reordering the $\left\{
\theta_{3},..,\theta_{n-6}\right\}  $ as before we obtain two even dimensional
Lie algebras, respectively isomorphic to $\frak{g}_{2m}^{13}$ and
$\frak{g}_{2m}^{26}.$

\item $\alpha_{i}^{4}=0,\;\forall i.$ A similar reordering of the $\theta_{i}$
gives two Lie algebras in odd dimension isomorphic to $\frak{g}_{2m+1}^{14}$
and $\frak{g}_{2m+1}^{27}.$
\end{enumerate}
\end{enumerate}

\item  Suppose now $\alpha_{1}^{3}=0.$

\begin{enumerate}
\item $\alpha_{1}^{4}=1,\;\alpha_{i}^{4}=0$ $\forall i\geq2.$ A linear change
allows to annihilate $\beta_{2}.$

\begin{enumerate}
\item $\alpha_{2}^{3}=1$ , and the remaining zero.\newline There are two
possible cases, depending on $\alpha_{23}^{1}:$ if it is nonzero we obtain an
odd dimensional algebra isomorphic to $\frak{g}_{2m+1}^{15}$ and if it is zero
an algebra in even dimension isomorphic to $\frak{g}_{2m}^{16}.$\newline 

\item $\alpha_{i}^{3}=0$ for any $i\geq2$ : we obtain an odd dimensional Lie
algebra isomorphic to $\frak{g}_{2m+1}^{17}.$
\end{enumerate}

\item $\alpha_{1}^{4}=0$

\begin{enumerate}
\item $\alpha_{2}^{3}\neq0$ and $\alpha_{i}^{3}=0,\forall i\geq3.$ With a
linear change we can suppose $\alpha_{2}^{4}=0.$\newline A-1) $\alpha_{3}%
^{4}=1$ and $\alpha_{i}^{4}=0$ for $i\geq4.$ A linear change allows to suppose
$\alpha_{3j}^{1}=0$ for all $j.$ If $\alpha_{2j}^{1}=0$ for all \ $j$ we
obtain the algebras $\frak{g}_{2m+1}^{18}$ and $\frak{g}_{2m+1}^{28}$. If not,
reorder $\left\{  \theta_{4},..,\theta_{n-6}\right\}  $ such that $\alpha
_{24}^{1}=1.$ We obtain the algebras $\frak{g}_{2m}^{9},\frak{g}_{2m}^{29}%
.$\newline A-2) $\alpha_{i}^{4}=0\;$ $\forall i.$ If $\alpha_{23}^{1}=0$ we
obtain the Lie algebras $\frak{g}_{2m}^{20},\frak{g}_{2m}^{30},$ and if
$\alpha_{23}^{1}\neq0$ we obtain the algebras $\frak{g}_{2m+1}^{21}$ and
$\frak{g}_{2m+1}^{31}.$

\item $\alpha_{i}^{3}=0\;\forall i.$\newline B-1) $\alpha_{2}^{4}=1,$%
\ $\alpha_{2j}^{1}=0.$ We obtain two Lie algebras in even dimension isomorphic
to $\frak{g}_{2m}^{22}$ and $\frak{g}_{2m}^{32}.$\newline B-2) $\alpha_{i}%
^{4}=0$ for $i\geq2$ : there is only one algebra in odd dimension which is
isomorphic to $\frak{g}_{2m+1}^{23}.$
\end{enumerate}
\end{enumerate}
\end{enumerate}
\end{proof}

There is only one remaining case, namely the corresponding to the $\left(
n-5\right)  $-filiform $2$-abelian Lie algebras with minimal dimension of its
derived algebra.

\begin{proposition}
Let $\frak{g}_{n}$ be an $2$-abelian algebra with $\dim\;C^{1}\frak{g}_{n}=4.$
Then $\frak{g}_{n}$ is isomorphic to one of the following laws :\newline 

\begin{enumerate}
\item $\frak{g}_{2m}^{33}\left(  m\geq4\right)  :$%
\[
\frak{g}_{0}^{2m}+\varphi_{1}+\psi_{1}^{3}+\psi_{2}^{4}+\sum
_{\substack{t=2\\m>4}}^{m-3}\psi_{2t-1,2t}^{1}%
\]

\item $\frak{g}_{2m}^{34}\left(  m\geq4\right)  :$%
\[
\frak{g}_{0}^{2m}+\varphi_{1}+\varphi_{2}+\psi_{1}^{3}+\psi_{2}^{4}%
+\sum_{\substack{t=2\\m>4}}^{m-3}\psi_{2t-1,2t}^{1}%
\]

\item $\frak{g}_{2m+1}^{35}\left(  m\geq4\right)  :$%
\[
\frak{g}_{0}^{2m+1}+\varphi_{1}+\psi_{1}^{3}+\psi_{2}^{4}+\psi_{13}^{1}%
+\sum_{\substack{t=2\\m>4}}^{m-3}\psi_{2t,2t+1}^{1}%
\]

\item $\frak{g}_{2m+1}^{36}\left(  m\geq4\right)  :$%
\[
\frak{g}_{0}^{2m+1}+\varphi_{1}+\varphi_{2}+\psi_{1}^{3}+\psi_{2}^{4}%
+\psi_{13}^{1}+\sum_{\substack{t=2\\m>4}}^{m-3}\psi_{2t,2t+1}^{1}%
\]

\item $\frak{g}_{2m+1}^{37}\left(  m\geq3\right)  :$%
\[
\frak{g}_{0}^{2m+1}+\varphi_{1}+\psi_{1}^{3}+\psi_{1}^{4}+\sum
_{\substack{t=1\\m>3}}^{m-3}\psi_{2t,2t+1}^{1}%
\]

\item $\frak{g}_{2m}^{38}\left(  m\geq4\right)  :$%
\[
\frak{g}_{0}^{2m}+\varphi_{1}+\psi_{1}^{3}+\sum_{t=1}^{m-3}\psi_{2t-1,2t}^{1}%
\]

\item $\frak{g}_{2m}^{39}\left(  m\geq4\right)  :$%
\[
\frak{g}_{0}^{2m}+\varphi_{1}+\varphi_{2}+\psi_{1}^{3}+\sum_{t=1}^{m-3}%
\psi_{2t-1,2t}^{1}%
\]

\item $\frak{g}_{2m+1}^{40}\left(  m\geq3\right)  :$%
\[
\frak{g}_{0}^{2m+1}+\varphi_{1}+\psi_{1}^{3}+\sum_{t=1}^{m-3}\psi
_{2t,2t+1}^{1}%
\]

\item $\frak{g}_{2m+1}^{41}\left(  m\geq3\right)  :$%
\[
\frak{g}_{0}^{2m+1}+\varphi_{1}+\varphi_{2}+\psi_{1}^{3}+\sum
_{\substack{t=1\\m>3}}^{m-3}\psi_{2t,2t+1}^{1}%
\]

\item $\frak{g}_{2m+1}^{42}\left(  m\geq3\right)  :$%
\[
\frak{g}_{0}^{2m+1}+\varphi_{1}+\psi_{1}^{4}+\sum_{\substack{t=1\\m>3}%
}^{m-3}\psi_{2t,2t+1}^{1}%
\]

\item $\frak{g}_{2m+1}^{43}\left(  m\geq3\right)  :$%
\[
\frak{g}_{0}^{2m+1}+\varphi_{1}+\varphi_{2}+\psi_{1}^{4}+\sum
_{\substack{t=1\\m>3}}^{m-3}\psi_{2t,2t+1}^{1}%
\]

\item $\frak{g}_{2m}^{44}\left(  m\geq3\right)  :$%
\[
\frak{g}_{0}^{2m}+\varphi_{1}+\sum_{\substack{t=1\\m>3}}^{m-3}\psi
_{2t-1,2t}^{1}%
\]

\item $\frak{g}_{2m}^{45}\left(  m\geq3\right)  :$%
\[
\frak{g}_{0}^{2m}+\varphi_{1}+\varphi_{2}+\sum_{\substack{t=1\\m>3}}^{m-3}%
\psi_{2t-1,2t}^{1}%
\]
\end{enumerate}
\end{proposition}

\begin{proof}
The starting assumptions for this case are
\[
\beta_{1}\neq0\;\text{\ and \ }\beta_{3,k}=0,k\geq1
\]
In particular the Jacobi condition forces $\alpha_{i}^{2}=0$ $\forall i.$

\begin{enumerate}
\item  Suppose $\beta_{1}$ and $\alpha_{1}^{3}\neq0$ $\left(  \text{so }%
\alpha_{j}^{3}=0\text{ for }i\geq2\text{ }\right)  .$\newline We observe that
if there exists an $\alpha_{ij}^{1}\neq0$ then a linear change of basis allows
to suppose $\alpha_{i}^{4}=\alpha_{j}^{4}=0.$ So we have the conditions
\begin{equation}
d_{i}a_{ij}=d_{j}a_{ij}=0,\;1\leq i,j
\end{equation}

\begin{enumerate}
\item $\exists\,\alpha_{i}^{4}\neq0$ with $i\geq2$ .We can suppose $\alpha
_{2}^{4}=1$ $\left(  \text{so }\alpha_{2j}^{1}=0\text{ by }\left(  1\right)
\right)  $ and $\alpha_{i}^{4}=0,\forall i\geq2.$

\begin{enumerate}
\item  If $\alpha_{1j}^{1}=0$ for all $j$ we obtain two even dimensional
algebras isomorphic to $\frak{g}_{2m}^{33}$ and $\frak{g}_{2m}^{34}.$

\item  If $\alpha_{1j}^{1}\neq0$ for an index $j$ we can suppose $\alpha
_{13}^{1}=1.$ We obtain two odd dimensional algebras isomorphic respectively
to $\frak{g}_{2m+1}^{35}$ and $\frak{g}_{2m+1}^{36}.$
\end{enumerate}

\item $\alpha_{i}^{4}$ $=0,\;\forall i\geq2$

\begin{enumerate}
\item  If $\alpha_{1}^{4}\neq0,$ then $\alpha_{1j}^{1}=0$ by $\left(
1\right)  .$ Reordering $\{\theta_{2},..,\theta_{n-6}\}$ we obtain an algebra
isomorphic to $\frak{g}_{2m+1}^{37}.$

\item  If $\alpha_{1}^{4}=0$ and $\alpha_{12}^{1}\neq0$ we obtain two algebras
isomorphic to $\frak{g}_{2m}^{38}$ and $\frak{g}_{2m}^{39}.$

\item  If $\alpha_{1}^{4}=\alpha_{1j}^{1}=0,\;\forall j$ we obtain two
algebras in odd dimension isomorphic to $\frak{g}_{2m+1}^{40}$ and
$\frak{g}_{2m+1}^{41}$.
\end{enumerate}
\end{enumerate}

\item  Suppose $\beta_{1}\neq0$ and $\alpha_{i}^{3}=0,\;\forall i.$%
\ Additionally we can suppose $\beta_{3}=0.$

\begin{enumerate}
\item  If $\alpha_{i}^{4}\neq0$ for an index $i\geq1$ let $\alpha_{1}^{4}=1$
and $\alpha_{i}^{4}=0,\;\forall i\geq2$ and $\alpha_{1j}^{1}=0$ by $\left(
1\right)  .$ We obtain two algebras isomorphic $\frak{g}_{2m+1}^{42}$ and
$\frak{g}_{2m+1}^{43}.$

\item  If $d_{i}=0$ $\forall i$ we obtain two even dimensional algebras
respectively isomorphic to $\frak{g}_{2m}^{44}$ and $\frak{g}_{2m}^{45}.$
\end{enumerate}
\end{enumerate}
\end{proof}

\begin{remark}
Observe that the algebras $\frak{g}_{2m+1}^{42}$ and $\frak{g}_{2m+1}^{43}$
are central extensions of the five dimensional filiform algebras $\frak{l}%
_{5}^{1}$ and $\frak{l}_{5}^{2}.$
\end{remark}

\begin{corollary}
Any nonsplit $\left(  n-5\right)  $-filiform $2$-abelian Lie algebra is
isomorphic to one of the laws $\frak{g}^{i},\;i\in\{1,..,45\}.$
\end{corollary}

\begin{remark}
As we have seen that a $\left(  n-5\right)  $-filiform Lie algebra is either
$1$- or $2$-abelian, the global classification follows from determining the
isomorphism classes of the $1$-abelian ones.
\end{remark}

\subsection{Characteristically nilpotent $\left(  n-5\right)  $-filiform Lie algebras}

The first example of a nilpotent Lie algebra all whose derivations are
nilpotent was given by Dixmier and Lister in 1957 [DL], as an answer to a
question formulated by Jacobson [Ja] two years earlier. This new class of Lie
algebras was soon recognized to be very important, and called
characteristically nilpotent, as they verify a certain sequence for
derivations which is a kind of generalization of the central descending
sequence for nilpotent Lie algebras ( [DL], [LT]).

\begin{definition}
A Lie algebra $\frak{g}$ is called characteristically nilpotent if the Lie
algebra of derivations $Der\left(  \frak{g}\right)  $ is nilpotent.
\end{definition}

\begin{remark}
It is easily seen that the original definition given by Dixmier and Lister 
is equivalent to the given above [LT].
\end{remark}

\begin{remark}
It is trivial to verify that there do not exist characteristically nilpotent,
$\left(  n-p\right)  $-filiform Lie algebras for indexes $p=1,2$. For $p=3,4$,
it has been shown that these algebras have rank $r\geq1$ [AC1], and that
almost any of these laws is the nilradical of a solvable, rigid law.
\end{remark}

\begin{lemma}
Let $\frak{g}$ be a $\left(  n-5\right)  $-filiform, $1$-abelian Lie algebra.
Then $rank\left(  \frak{g}\right)  \geq1$.
\end{lemma}

\begin{proof}
If $\dim C^{1}\frak{g}=4$, the assertion follows immediately from the linear
system $\left(  S\right)  $ associated to the algebra, as this system admits
nontrivial solutions. If $\dim C^{1}\frak{g}=5$, the only case for which the
system could have zero solution is $\beta_{2}=\beta_{3,1}=1$, and the distinct
values of $\left(  \alpha_{i}^{2},\alpha_{i}^{3},\alpha_{i}^{4}\right)  $. For
any of these starting conditions it is routine to prove the existence of a
nonzero semisimple derivation.
\end{proof}

\begin{lemma}
Let $\frak{g}\in\frak{H}_{2}$. If $\alpha_{ij}^{1}\neq0$ for $1\leq i,j\leq
n-6$ such that
\[
\beta_{1,k}=\beta_{3,k}=\alpha_{k}^{t}=0,\;k=i,j,\;t=2,3,4
\]
then $rank\left(  \frak{g}\right)  \geq1.$
\end{lemma}

\begin{proof}
Consider the endomorphism defined by
\[
d\left(  Y_{i}\right)  =Y_{i},\;d\left(  Y_{j}\right)  =-Y_{j}%
\]
and zero over the undefined images, where $\left(  X_{1},..,X_{6}%
,Y_{1},..,Y_{n-6}\right)  $ is the dual basis of $\left(  \omega_{1}%
,..,\omega_{6},\theta_{1},..,\theta_{n-6}\right)  $. Clearly $d$ is a nonzero
semisimple derivation of $\frak{g}$.
\end{proof}

\begin{proposition}
A $\left(  n-5\right)  $-filiform Lie algebra $\frak{g}_{n}$ is
characteristically nilpotent if and only if it is isomorphic to one of the
following laws:%
\begin{align*}
&  \frak{g}_{7}^{11},\;\frak{g}_{9}^{15},\;\;\frak{g}_{7}^{17},\;\frak{g}%
_{7}^{24,\alpha}\,\left(  \alpha\neq0\right)  ,\;\frak{g}_{7}^{27}%
,\;\frak{g}_{9}^{36},\;\;\frak{g}_{7}^{37},\;\frak{g}_{7}^{41}\\
&  \frak{g}_{8}^{12},\;\frak{g}_{8}^{16},\;\frak{g}_{8}^{25,\alpha}\,\left(
\alpha\neq0\right)  ,\;\frak{g}_{8}^{26},\;\frak{g}_{8}^{34},\;\frak{g}%
_{8}^{39}%
\end{align*}
\end{proposition}

\begin{corollary}
There are characteristically nilpotent Lie algebras $\frak{g}_{n}$ with
nilpotence index $5$ for the dimensions $n=7,8,9,14,15,16,17,18$ \ and
\ $n\geq21.$
\end{corollary}

\begin{proof}
As the sum of characteristically nilpotent algebras is characteristically
nilpotent [T1], the assertion follows from the previous proposition.
\end{proof}

\begin{remark}
In fact, for any $n\geq7$ there exist characteristically nilpotent Lie
algebras of nilindex $5$. However, the algebras to be added are not
$p$-filiform any more [AC2].
\end{remark}

\subsection{Nilradicals of rigid algebras as factors of $k$-abelian Lie algebras}

The second application of $k$-abelian Lie algebras is of interest for the
theory of rigid Lie algebras. In this paragraph we prove the existence, by
giving a family for dimensions $2m+2\;\left(  m\geq4\right)  $, of $\left(
m-1\right)  $-abelian Lie algebras $\frak{g}$ of characteristic sequence
$\left(  2m-1,2,1\right)  $ all whose factor algebras $\frac{\frak{g}}%
{C^{k}\frak{g}}$ $\left(  k\geq m\right)  $ are isomorphic to the nilradical
of a solvable rigid law.

For $m\geq4$ let $\frak{g}_{m}$ be the Lie algebra whose Cartan-Maurer
equations are
\begin{align*}
d\omega_{1}  &  =d\omega_{2}=0\\
d\omega_{j}  &  =\omega_{1}\wedge\omega_{j-1},\;3\leq j\leq2m-1\\
d\omega_{2m}  &  =\omega_{1}\wedge\omega_{2m-1}+\sum_{j=2}^{m}\left(
-1\right)  ^{j}\omega_{j}\wedge\omega_{2m+1-j}\\
d\omega_{2m+1}  &  =\omega_{2}\wedge\omega_{3}\\
d\omega_{2m+2}  &  =\omega_{1}\wedge\omega_{2m+1}+\omega_{2}\wedge\omega_{4}%
\end{align*}
It is elementary to verify that this algebra has characteristic sequence
$\left(  2m-1,2,1\right)  $. Moreover, it is $\left(  m-1\right)  $-abelian,
for the exterior product $\omega_{m}\wedge\omega_{m+1}$ proves that $\left[
C^{m-2}\frak{g}_{m},C^{m-2}\frak{g}_{m}\right]  \neq0$ and $\left[
C^{m-1}\frak{g}_{m},C^{m-1}\frak{g}_{m}\right]  =0.$

\begin{notation}
The dual basis of $\left(  \omega_{1},..,\omega_{m2+2}\right)  $ will be
denoted as $\left(  X_{1},..,X_{2m+2}\right)  $.
\end{notation}

\begin{lemma}
For any $4\leq m\leq k\leq2m-2$ the factor algebra $\frac{\frak{g}_{m}}%
{C^{k}\frak{g}_{m}}$ has equations
\begin{align*}
d\overset{-}{\omega}_{1} &  =d\overset{-}{\omega}_{2}=0\\
d\overset{-}{\omega}_{j} &  =\overset{-}{\omega}_{1}\wedge\overset{-}{\omega
}_{j-1},\;3\leq j\leq k+1\\
d\overset{-}{\omega}_{2m+1} &  =\overset{-}{\omega}_{2}\wedge\overset
{-}{\omega}_{3}\\
d\overset{-}{\omega}_{2m+2} &  =\overset{-}{\omega}_{1}\wedge\overset
{-}{\omega}_{2m+1}+\overset{-}{\omega}_{2}\wedge\overset{-}{\omega}_{4}%
\end{align*}
where $\overset{-}{\omega_{j}}=\omega_{j}\;\operatorname{mod}C^{k}\frak{g}%
_{m}$. Moreover, this algebra is $1$-abelian of characteristic sequence
$\left(  k,2,1\right)  $.
\end{lemma}

The proof is trivial.

\begin{proposition}
For any $4\leq m\leq k$ the algebra $\frac{\frak{g}_{m}}{C^{k}\frak{g}_{m}}$
is isomorphic to the nilradical of a solvable, rigid Lie algebra
$\frak{r}_{m,k}$.
\end{proposition}

\begin{proof}
Let $\frak{r}_{m,k}=$ $\frac{\frak{g}_{m}}{C^{k}\frak{g}_{m}}\oplus
\frak{t}_{k}$ be the semidirect product of $\frac{\frak{g}_{m}}{C^{k}%
\frak{g}_{m}}$ by the torus $\frak{t}_{k}$ defined by its weights :
\[%
\begin{tabular}
[c]{c}%
$\lambda_{1},\;\lambda_{2}+\left(  k-1\right)  \lambda_{1},\;\lambda
_{j}=\lambda_{2}+\left(  k-3+j\right)  \lambda_{1}\;\left(  3\leq j\leq
k+1\right)  $\\
$\;\lambda_{2m+1}=2\lambda_{2}+\left(  2k-1\right)  \lambda_{1},\;\lambda
_{2m}=2\lambda_{2}+2k\lambda_{1}$%
\end{tabular}
\]
over the basis $\left\{  \overset{-}{X}_{1},..,\overset{-}{X}_{k+1}%
,\overset{-}{X}_{2m+1},\overset{-}{X}_{2m+2}\right\}  $ dual to $\left\{
\overset{-}{\omega}_{1},..,\overset{-}{\omega}_{k+1},\overset{-}{\omega
}_{2m+1},\overset{-}{\omega}_{2m+2}\right\}  $. Then the law is given by%
\[%
\begin{tabular}
[c]{l}%
$\left[  V_{1},\overset{-}{X}_{1}\right]  =\overset{-}{X}_{1},\;\left[
V_{1},\overset{-}{X}_{j}\right]  =\left(  k-1\right)  \overset{-}{X}%
_{j}\left(  3\leq j\leq k+1\right)  ,$\\
$\left[  V_{1},\overset{-}{X}_{2m+1}\right]  =\left(  2k-1\right)  \overset
{-}{X}_{2m+1},\;\left[  V_{1},\overset{-}{X}_{2m+2}\right]  =2k\overset{-}%
{X}_{2m+2}$\\
$\left[  V_{2},\overset{-}{X}_{j}\right]  =\overset{-}{X}_{j}\left(  2\leq
j\leq k+1\right)  ,\;\left[  V_{2},\overset{-}{X}_{j}\right]  =2\overset{-}%
{X}_{j}\left(  j=2m+1,2m+2\right)  $\\
$\left[  \overset{-}{X}_{1},\overset{-}{X}_{j}\right]  =\overset{-}{X}%
_{j+1}\left(  2\leq j\leq k\right)  ,\;\left[  \overset{-}{X}_{2},\overset
{-}{X}_{3}\right]  =a\overset{-}{X}_{2m+1},\;\left[  \overset{-}{X}%
_{2},\overset{-}{X}_{4}\right]  =b\overset{-}{X}_{2m+1}$%
\end{tabular}
\]

Thus the only nonzero brackets not involving the vector $\overset{-}{X}_{1}$
are
\[
\left[  \overset{-}{X}_{2},\overset{-}{X}_{3}\right]  =a\overset{-}{X}%
_{2m+1},\;\left[  \overset{-}{X}_{2},\overset{-}{X}_{4}\right]  =b\overset
{-}{X}_{2m+1}%
\]
Now Jacobi implies $a=b$, and by a change of basis $a=1$. Thus the law is
rigid, and $\frak{t}_{k}$ is a maximal torus of derivations of $\frac
{\frak{g}_{m}}{C^{k}\frak{g}_{m}}$.
\end{proof}

\begin{corollary}
For any $4\leq m\leq k$ the factor algebra
\[
\frac{\left(  \frac{\frak{g}_{m}}{C^{k}\frak{g}_{m}}\right)  }{\left\langle
\overset{-}{X}_{2m+2}\right\rangle }%
\]
is $1$-abelian of characteristic sequence $\left(  k,1,1\right)  $ and
isomorphic to the nilradical of a solvable rigid law $\frak{s}_{m,k}$.
Moreover
\[
\frak{s}_{m,k}\simeq\frac{\frak{r}_{m,k}}{\left\langle \overset{-}{X}%
_{2m+2}\right\rangle }%
\]
\end{corollary}

\section{Other 2-abelian nilpotent Lie algebras}

Let $E_{6}$ be the simple exceptional Lie algebra of dimension $78$. Let
$\Phi$ be a root system respect to a Cartan subalgebra $\frak{h}$ and
$\Delta=\left\{  \alpha_{1},..,\alpha_{6}\right\}  $ a basis of fundamental
roots. Recall that the standard Borel subalgebra of $E_{6}$ is given by
\[
B\left(  \Delta\right)  =\frak{h}+\sum_{\alpha\in\Phi^{+}}L_{\alpha}%
\]
where $L_{\alpha}$ is the weight space associated to the root $\alpha$. Recall
also that any parabolic subalgebra $\frak{p}$ is determined, up to isomorphism, by a
subsystem $\Delta_{1}\subset\Delta$ such that $P$ is conjugated to the
subalgebra
\[
P\left(  \Delta_{1}\right)  =\frak{h}+\sum_{\alpha\in\Phi_{1}\cup\Phi_{2}^{+}%
}L_{\alpha}%
\]
where $\Phi_{1}$ is the set of roots expressed in terms of $\Delta
\backslash\Delta_{1}$ and $\Phi_{2}^{+}=\Phi^{+}\cap\left(  \Phi\backslash
\Phi_{1}\right)  $. It is elementary to see that the nilradical is
\[
\frak{n}\left(  \Delta_{1}\right)  =\sum_{\alpha\in\Phi_{2}^{+}}L_{\alpha}%
\]
and called $\left(  E_{6},\Delta_{1}\right)  $-nilalgebra.

Let
\[
\mathcal{L}=\left\{
\begin{array}
[c]{c}%
\left\{  \alpha_{1},\alpha_{4}\right\}  ,\left\{  \alpha_{4},\alpha
_{6}\right\}  ,\left\{  \alpha_{3},\alpha_{5}\right\}  ,\left\{  \alpha
_{3},\alpha_{4}\right\}  ,\left\{  \alpha_{4},\alpha_{5}\right\}  ,\left\{
\alpha_{2},\alpha_{3}\right\}  ,\\
\left\{  \alpha_{2},\alpha_{5}\right\}  ,\left\{  \alpha_{2},\alpha
_{4}\right\}  ,\left\{  \alpha_{1},\alpha_{2},\alpha_{3}\right\}  ,\left\{
\alpha_{2},\alpha_{5},\alpha_{6}\right\}  ,\\
\left\{  \alpha_{1},\alpha_{2},\alpha_{5}\right\}  ,\left\{  \alpha_{2}%
,\alpha_{3},\alpha_{6}\right\}  ,\left\{  \alpha_{1},\alpha_{4},\alpha
_{6}\right\}  ,\\
\left\{  \alpha_{1},\alpha_{2},\alpha_{6}\right\}  ,\left\{  \alpha_{1}%
,\alpha_{3},\alpha_{5}\right\}  ,\left\{  \alpha_{3},\alpha_{5},\alpha
_{6}\right\}
\end{array}
\right\}
\]

\bigskip As known, the maximal root of $E_{6}$ is
\[
\delta=\alpha_{1}+2\alpha_{2}+2\alpha_{3}+3\alpha_{4}+2\alpha_{5}+\alpha
_{6}=\sum_{i=1}^{6}k_{a_{i}}\alpha_{i}%
\]
We define the $\Delta_{1}$-height of $\delta$ as
\[
h_{\Delta_{1}}\left(  \delta\right)  =\sum_{\alpha_{i}\in\Delta_{1}}k_{a_{i}}%
\]
and the subsets
\[
\Delta_{1}\left(  k\right)  =\left\{  \alpha\in\Phi_{2}^{+}\;|\;h_{\Delta_{1}%
}\left(  \alpha\right)  =k\right\}
\]

\begin{proposition}
For any $\Delta_{1}\in\mathcal{L}$ the $\left(  E_{6},\Delta_{1}\right)
$-nilalgebra $\frak{n}\left(  \Delta_{1}\right)  $ is $2$-abelian.
\end{proposition}

\begin{proof}
As known, for the ideals $C^{k}\frak{n}$ of the descending central sequence we
have
\[
C^{k}\frak{n}=\sum_{\substack{\alpha\in\Delta_{1}\left(  j\right)  \\j\geq
k+1}}L_{\alpha}%
\]
Thus, if the derived subalgebra is not abelian, it suffices to show the
existence of two roots $\alpha,\beta\in\Delta_{1}\left(  2\right)  $ such that
$\alpha+\beta\in\Phi_{2}^{+}$ and that for any two roots $\gamma
,\varepsilon\in\Delta_{1}\left(  3\right)  $ we have $\gamma+\varepsilon
\notin\Phi_{2}^{+}.$ Moreover, let $\delta_{1}=\sum_{i=1}^{6}\alpha_{i}\in\Phi$

\begin{enumerate}
\item $\Delta_{1}=\left\{  \alpha_{1},\alpha_{4}\right\}  $ : take
$\alpha=\delta_{1}-\alpha_{5}-\alpha_{6},\;\beta=\delta-\alpha_{1}-\alpha
_{2}-\alpha_{4};\;\alpha+\beta=\delta$

\item $\Delta_{1}=\left\{  \alpha_{3},\alpha_{5}\right\}  $ : $\alpha
=\delta_{1},\;\beta=\delta_{1}-\alpha_{1}-\alpha_{6}+\alpha_{4};\;\alpha
+\beta=\delta$

\item $\Delta_{1}=\left\{  \alpha_{4},\alpha_{5}\right\}  $ : $\alpha
=\delta_{1}-\alpha_{1}-\alpha_{2}-\alpha_{6},\;\beta=\delta_{1};\;\alpha
+\beta=\delta-\alpha_{2}-\alpha_{4}$

\item $\Delta_{1}=\left\{  \alpha_{2},\alpha_{3}\right\}  $ : $\alpha
=\delta_{1},\;\beta=\delta_{1}-\alpha_{1}-\alpha_{2}-\alpha_{6};\;\alpha
+\beta=\delta-\alpha_{2}-\alpha_{4}$

\item $\Delta_{1}=\left\{  \alpha_{2},\alpha_{4}\right\}  $ : $\alpha
=\delta_{1},\;\beta=\delta_{1}-\alpha_{1}-\alpha_{6}+\alpha_{4};\;\alpha
+\beta=\delta$

\item $\Delta_{1}=\left\{  \alpha_{1},\alpha_{2},\alpha_{3}\right\}  $ :
$\alpha=\delta-\delta_{1},\;\beta=\delta_{1}-\alpha_{2};\;\alpha+\beta
=\delta-\alpha_{2}$

\item $\Delta_{1}=\left\{  \alpha_{1},\alpha_{2},\alpha_{5}\right\}  $ :
$\alpha=\delta-\delta_{1},\;\beta=\delta_{1}-\alpha_{2};\;\alpha+\beta
=\delta-\alpha_{2}$

\item $\Delta_{1}=\left\{  \alpha_{1},\alpha_{3},\alpha_{6}\right\}  $ :
$\alpha=\delta_{1}-\alpha_{6},\;\beta=\delta_{1}-\alpha_{1}-\alpha
_{2};\;\alpha+\beta=\delta-\alpha_{2}-\alpha_{4}$

\item $\Delta_{1}=\left\{  \alpha_{1},\alpha_{4},\alpha_{6}\right\}  $ :
$\alpha=\delta_{1}-\alpha_{6},\;\beta=\delta_{1}-\alpha_{1}-\alpha
_{2};\;\alpha+\beta=\delta-\alpha_{2}-\alpha_{4}$

\item $\Delta_{1}=\left\{  \alpha_{1},\alpha_{2},\alpha_{6}\right\}  $ :
$\alpha=\delta_{1}-\alpha_{1},\;\beta=\delta_{1}+\alpha_{4};\;\alpha
+\beta=\delta$

\item $\Delta_{1}=\left\{  \alpha_{1},\alpha_{3},\alpha_{5}\right\}  $ :
$\alpha=\delta_{1}-\alpha_{2}-\alpha_{6},\;\beta=\delta_{1}+\alpha
_{4};\;\alpha+\beta=\delta-\alpha_{2}$
\end{enumerate}

Let $X_{\alpha},X_{\beta}$ be generators of the weight spaces $L_{\alpha}$ and
$L_{\beta}$ : then the preceding relations show that
\[
\left[  X_{\alpha},X_{\beta}\right]  =X_{\alpha+\beta}\neq0
\]
proving that the derived subalgebra is not abelian.\newline Finally, it is
trivial to see that for any subset $\Delta_{1}$ listed above we have
\[
\left[  C^{2}\frak{n}\left(  \Delta_{1}\right)  ,C^{2}\frak{n}\left(
\Delta_{1}\right)  \right]  =0
\]
\end{proof}


\begin{thebibliography}{99}
\bibitem[ AG1]{}J. M. Ancochea, M. Goze. \textit{On the nonrationality of
rigid Lie algebras, }Proc. Amer. Math. Soc. 129(9), 2611-2618 (1999).

\bibitem[ AG2] {}M. Goze, J. M. Ancochea. \textit{Alg\`{e}bres de Lie rigides,
}Proceedings A 88(4), 397-415 (1985).

\bibitem[ AG3] {}J. M. Ancochea, M. Goze. \textit{Sur la classification des
alg\`{e}bres de Lie nilpotentes de dimension 7, }C. R. Acad. Sci. Paris 302
(1986), 611-613.

\bibitem[ AC1] {}J. M. Ancochea, O. R. Campoamor. \textit{On Lie algebras
whose nilradical is (n-p)-filiform, }Comm. Algebra (29), 2001, 427-450.

\bibitem[ AC2] {}J. M. Ancochea, O. R. Campoamor\textit{. Characteristically nilpotent Lie algebras, } Cont. Matematicas 3 (2000).

\bibitem[ Bra] {}F. Bratzlavsky. \textit{Sur les alg\`{e}bres admettant un
tores des d\'{e}rivations donn\'{e}, }Journal of Algebra (30), 1974, 304-316.

\bibitem[ Ca] {}R. Carles. \textit{Sur les alg\`{e}bres
caract\'{e}ristiquement nilpotentes. }Publ. Univ. Poitiers 1984.

\bibitem[ DiL] {}J. Dixmier, W. G. Lister. \textit{Derivations of nilpotent
Lie algebras, }Proc. Amer. Math. Soc. 8 (1957), 155-157.

\bibitem[ Di] {}J. Dixmier. \textit{Cohomologie des alg\`{e}bres de Lie
nilpotentes, }Acta Sci. Math. Szeged 16 (1955), 246-250.

\bibitem[ Dy] {}J. L. Dyer. \textit{A nilpotent Lie algebra with nilpotent
automorphism group, }Bull. Amer. Math. Soc. 76 (1970), 52-56.

\bibitem[ Ei] {}S. Eilenberg. \textit{Extensions of general algebras, }Ann.
Soc. Polon. Math. 21 (1948), 125-134.

\bibitem[ Fa] {}G. Favre. \textit{Un alg\`{e}bre de Lie
caract\'{e}ristiquement nilpotente en dim 7, }C. R. Acad. Sci. Paris 274
(1972), 1338-1339.

\bibitem[ GaTm] {}L. Yu. Galitski, D. A. Timashev. \textit{On classifications
of metabelian Lie algebras, }J. of Lie theory 9 (1999), 125-156.

\bibitem[ GoGK] {}J. R. G\'{o}mez, M. Goze, Yu. B. Khakimdjanov. \textit{On
the k-abelian filiform Lie algebras. }Comm. Algebra 25(2), 1997, 431-450.

\bibitem[ GK1] {}M. Goze, Yu. B. Khakimdjanov. \textit{Nilpotent Lie algebras,
}Kluwer Ac. Press 1996.

\bibitem[ G] {}M. Goze. \textit{Perturbations sur une vari\'{e}t\'{e}
alg\'{e}brique: application \`{a} l'\'{e}tude de la vari\'{e}t\'{e} des lois
d'alg\`{e}bre de Lie sur ${C}^{n}$, } C. R. A. Sci. Paris 295 (1982), 583-586.

\bibitem[ Ja] {}N. Jacobson. \textit{A note on automorphisms and derivations
of Lie algebras, }\ Proc.Amer. Math. Soc. 6 (1955), 281-283.

\bibitem[ LeLu1] {}\ G. F. Leger, E. Luks. \textit{On a duality for metabelian
Lie algebras, }J. of Algebra 21 (1972), 266-270.

\bibitem[ LeL2] {}G. F. Leger, E. Luks. \textit{On derivations and holomorphs
of nilpotent Lie algebras, }Nagoya Math. J. 44 (1971), 39-50.

\bibitem[ LeT] {}G. F. Leger, S. T\^{o}g\^{o}. \textit{Characteristically
nilpotent Lie algebras, }Duke Math. J. 26 (1959), 623-628.

\bibitem[ L] {}E. Luks. \textit{What is the typical nilpotent Lie algebra, }in
Computers in nonassociative rings and algebras, Acad. Press 1977.

\bibitem[ Mo] {}V. V. Morozov. \textit{Classification des alg\`{e}bres de Lie
nilpotentes de dimension 6, }Izv. Vyssh. Ucheb. Zar. 4 (1958), 161-171.

\bibitem[ Mos] {}G. D. Mostow. \textit{Fully reductible subgroups of algebraic
groups, }Amer. J. Math. 68 (1956), 200-221.

\bibitem[ Mu] {}F. J. Murray. \textit{Perturbation theory and Lie algebras,
}J. Math. Phys. 3 (1962), 451-468.

\bibitem[ NiR] {}A. Nijenhuis, R. W. Richardson. \textit{Deformations of Lie
algebra structures, }J. Math. Mech. 17 (1967), 89-105.

\bibitem[ Ri] {}D. S. Rim. \textit{Deformations of transitive Lie algebras,
}Ann. of Math. (1966), 339-357.

\bibitem[ Sch] {}E. Schenkmann. \textit{A theory of subinvariant Lie algebras,
}Amer. J. Math. 73 (1951), 453-474.

\bibitem[ T1] {}S. T\^{o}g\^{o}. \textit{Outer derivations od Lie algebras,
}Trans. Amer. Math. Soc. 128 (1967), 264-276.

\bibitem[ T2] {}S. T\^{o}g\^{o}. \textit{On the derivation algebras of Lie
algebras, }Cand. J. Math. 13(2), 201-216 (1961).

\bibitem[ Ve1] {}M. Vergne. \textit{Vari\'{e}t\'{e} des alg\`{e}bres de Lie
nilpotentes, }These 3$^{eme}$ \textit{cycle, Paris 1966.}

\bibitem[ Ve2] {}M. Vergne. \textit{Cohomologie des alg\`{e}bres de Lie
nilpotentes. Applications a l'\'{e}tude de la vari\'{e}t\'{e} des alg\`{e}bres
de Lie nilpotentes, }Bull. Soc. Math. France 98 (1970), 81-116.

\bibitem[ Ya1] {}S. Yamaguchi. \textit{Derivations and affine structures of
some nilpotent Lie algebras, }Mem. Fac. Sci. Kyushu Univ. Ser. A 34 (1980), 151-170.

\bibitem[ Ya2] {}S. Yamaguchi. \textit{On some classes of nilpotent Lie
algebras and their automorphism group, }Mem. Fac. Sci. Kyushu Univ. Ser. A 34
(1981), 241-251.

\bibitem[ Zh] {}L. Zhu \textit{The construction of some solvable complete Lie
algebras, }J. Nanjing Univ. Math. Biquat. 15(1), 34-40 (1998).
\end{thebibliography}
\end{document}